\documentclass[reqno, 12pt]{amsart}

\usepackage{graphics}
\usepackage[pdftex]{graphicx}           

\usepackage{amssymb}
\usepackage{latexsym}
\usepackage{epsfig}
\usepackage{epsf}
\usepackage{float}
\usepackage{a4wide}
\usepackage{hyperref}

\usepackage{color}

\newtheorem{theorem}{Theorem}
\newtheorem{lemma}{Lemma}

\newtheorem{proposition}{Proposition}

\def\bN{{\mathbb N}}
\def\bP{{\mathbb P}}

\def\bE{{\mathbb E}}

\def\0{{\mathbf 0}}

\def\1{{\mathbf 1}}

\def\cH{{\mathcal H}}

\def\bchi{{\overline \chi}}
\def\beeta{{\overline \eta}}
\def\bsigma{{\overline \sigma}}
\def\bxi{{\overline \xi}}
\def\bS{{\bar S}}
\def\bN{{\bar N}}
\def\bW{{\bar W}}
\def\bEE{{\bar E}}
\def\bI{{\bar I}}
\def\bJ{{\bar J}}

\def\var{{\rm Var}}
\def\cov{{\rm Cov}}

\def\reff#1{(\ref{#1})}
\def\proofof #1{{\noindent \bf Proof of #1.}}
\def\endproof{$\square$ \vskip 2mm}

\begin{document}

\title[The variance of the shock ]{The variance of the shock in the HAD process}
\author{Cristian F. Coletti, Pablo A. Ferrari and Leandro P. R. Pimentel}
\address{Instituto de Matem\'atica\\
  Universidade de S\~{a}o Paulo\\
  Rua do Mat\~{a}o 1010 bloco D 05311-970 S\~{a}o Paulo SP Brazil\\}
\email{cristian@ime.usp.br, pablo@ime.usp.br, ordnael@ime.usp.br} \urladdr{}

\subjclass[2000]{}


\begin{abstract}
  We consider the Hammersley-Aldous-Diaconis (HAD) process with sinks and
  sources such that there is a microscopic shock at every time $t$; denote
  $Z(t)$ its position. We show that the mean and variance of $Z(t)$ are linear
  functions of $t$ and compute explicitely the respective constants in function
  of the left and right densities. Furthermore, we describe the dependence of
  $Z(t)$ on the initial configuration in the scale $\sqrt t$ and, as a
  corollary, prove a central limit theorem.
 \end{abstract}

\maketitle

\section{Introduction}

\begin{figure}[h]
\begin{center}
\includegraphics[width=0.5\textwidth]{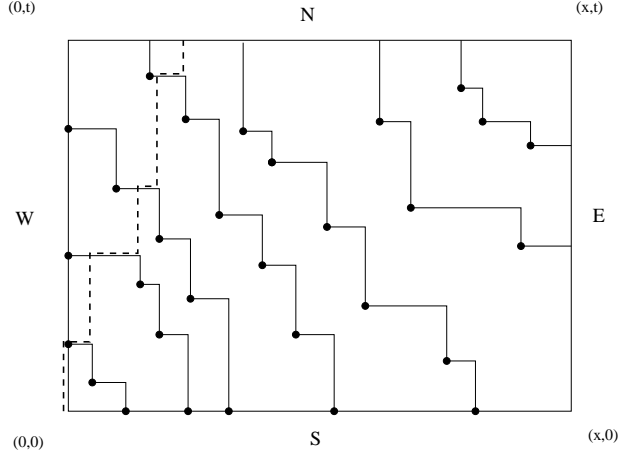}
\end{center}
\caption{The black points represent the Poisson points $S$ (sources), $W$
  (sinks) and $P$. The polygonal lines represent the trajectory of HAD particles
  while the trajectory of the second class particle is represented by the dashed
  polygonal line.}\label{HAD}
\end{figure}

Let $S$ and $W$ be one-dimensional Poisson processes and let $P$ be a two
dimensional Poisson process of rate 1. Assume they are homogeneous and mutually
independent.  The Hammersley-Aldous-Diaconis process \cite{h, ad}, or shortly
the HAD process, $\cH(S,W,P)=(H_s, \,s\in[0,t])$\ has been constructed by
Groeneboom \cite{g} in the square $[0,x]\times[0,t]$ as a deterministic function
of $S$, $W$ and $P$ as follows. The point configuration $H_s$ represents the
position of particles. At time zero the particles start at the points in $S$,
called the \emph{sources}. Then (see Figure \ref{HAD}), at the first $s>0$ such
that $(y,s)$ is in $P$ for some $y\in[0,x]$ or $s$ is in $W$, the closest
particle to the right of $y$ jumps to $y$ if $(y,s)$ is in $P$ or to $0$ if $s$
is in $W$ (points in $W$ are called \emph{sinks}). If there is no particle to
the right of $y$, then a new particle is added at $y$ at time $s$. The new
configuration does not move until the second $s$ such that $(y,s)$ is in $P$ for
some $y$ or $s$ is in $W$, when the second jump occurs, and so on until time
$t$. In other words, we define the process inductively as follows: $H_0=S$ and
for $s>0$,
\[
H_s =\left\{
\begin{array}{ll}
  H_{s-} &\hbox { if } s\notin W \cup \{s',\, (y,s)\in P \hbox{ for some }
  y'\in[0,x]\}\\
\Big\{ H_{s-} \setminus \{R(y,H_{s-})\}\Big\} \cup \{y\}&\hbox { if } (y,s)\in P \cup
\{(0,s'),\,s'\in W\}
\end{array}
\right.
\]
where $R(y,H) = \inf\{y'\in H\cup\{x\}\,:\, y'>y\}$ for $H\subset (0,x)$. Let
$N$ be the positions of particles at time $t$, and let $E$ be the set of times a
new particle enters the system from the right through the vertical axis
$\{x\}\times [0,t)$.

Let $\lambda\geq 0$ and $\rho\geq 0$. Assume that $S$ has intensity $\lambda$,
$W$ has intensity $\rho$ and $P$ has intensity $1$. Define the point process
$S'$ by removing the first sink of $S$ and consider the corresponding HAD
process $\cH(S',W,P)$. For time $t\geq 0$ the coupled processes $\cH(S,W,P)$ and
$\cH(S',W,P)$ will differ at one point denoted $Z(t)$ and called a second class
particle (Figure \ref{HAD}).

\paragraph{\bf Results} Our main results are the computation of the mean and
variance of $Z(t)$ and the dependence of fluctuations on the initial
configuration.

\begin{theorem}\label{2class-mean&var}
Assume that $\lambda\rho>1$. Then for all $t\geq 0$ we have
\begin{equation}\label{2class-mean}
\bE\big(Z(t)\big)\,=\,\frac{\rho}{\lambda}t\,,
\end{equation}
\begin{equation}\label{2class-var}
\var(Z(t))=  \frac{2\big(\rho-\frac{1}{\lambda}\big)}{\big(\lambda-\frac{1}{\rho}\big)^{2}}t\,.
\end{equation}
\end{theorem}

Let
$$
R(S)=\frac{1}{\lambda}\Big(\rho-\frac{1}{\lambda}\Big)\,\mbox{ and
}R(W)=\frac{1}{\rho}\Big(\rho-\frac{1}{\lambda}\Big)\,.
$$
In the next theorem we show that in the
scale $\sqrt t$ the position of $Z(t)$ is the same as the following random
variable, depending only on $S$ and $W$:
$$
N(t):=\big(\lambda-\frac{1}{\rho}\big)^{-1}\Big[3\Big(\rho-\frac{1}{\lambda}\Big)
t-\Big(\big|[0,tR(S)]\cap S\big|+\big|[0,tR(W)]\cap W\big|
\Big)\Big]\,,
$$
where $|A|$ denotes the number of points in the (finite) set $A$.

\begin{theorem}\label{dependence}
Assume that $\lambda\rho>1$. Then
\[
\lim_{t\to\infty}  \frac{\bE\Big[\big(Z(t)-N(t)\big)^2\Big]}{t} =0\,.
\]
and
\[
\lim_{t\to\infty}  \frac{Z(t)-\bE\big(Z(t)\big)}{\sqrt{\var \big(Z(t)\big)}}= \mathcal
N(0,1)\,\mbox{ in distribution}\,,
\]
where $\mathcal N(0,1)$ is a standard normal random variable.\end{theorem}



\paragraph{\bf Brief historical overview} The HAD process has a genealogical
relation with the famous Ulam's problem, which studies the limit behavior of the
longest increasing subsequence $L_n$ in a random permutation of the numbers
$1,\dots,n$. At the beginning of the seventies Hammersley \cite{h} proposed to
solve this problem as follows. Assume $(x_1,t_1),\dots,(x_n,t_n)$ are $n$
independent random points distributed uniformly in the rectangle
$[0,x]\times[0,t]$. These points specify a uniform random permutation $\pi$ by
setting the point with the $i$'th smallest $t$-coordinate has the $\pi(i)$
smallest $x$-coordinate. The length of the longest increasing subsequence of
$\pi$ equals the maximal number $M(x,t)$ of points on an up-right path from
$(0,0)$ to $(x,t)$ (last-passage percolation). The variable $M(x,t)$ has the so
called super-additivity property \cite{h+}, which was the key to show that
$\sqrt{n}L_n\to c$.  Hammersley conjectured that $c=2$ by presenting a
non-rigorous but quite reasonable argument based on hydrodynamical local
equilibrium for a related discrete-time particle system. Later Logan and Sheep
\cite{ls} and Vershik and Kerov \cite{kv} proved $c=2$ with combinatorial
methods.

In the middle of the nineties, Aldous and Diaconis \cite{ad}
constructed a continuous time version of the Hammersley process and
use Rost \cite{r} coupling ideas, developed in the exclusion process
context, to show that $c=2$. This is the reason we call this process
Hammersley-Aldous-Diaconis process. Groeneboom \cite{g} introduced
the construction in the quadrant we use here. In this setup $M(x,t)$
equals the number of particles in $[0,x]\times\{t\}$ plus the number
of particles in $\{0\}\times [0,t]$ (we are allowing an up-right
path to pick Poissonian points in the horizontal or vertical lines).
Thus, the trajectory of the $n$'th particle corresponds to the
boundary of the region $\{ M(x,t)\leq n-1\}$ (Figure 1).

The macroscopic hydrodynamical behavior of the HAD process is given by the
Burgers equation $\partial_t u+ \partial_x g(u) = 0$, with $g(u)=\frac{1}{u}$
\cite{s}. When the initial conditions are $u(r,0)=\lambda$ for $r>0$ and
$u(r,0)=\frac1\rho$ when $r\le 0$, the solution has a drastic changing when the
product $\lambda\rho$ goes from values smaller than $1$ to values greater than
$1$. For $\lambda\rho=1$, the system is stationary and the solutions are
constants in time. For $\lambda\rho<1$ it develops a rarefaction front while for
$\lambda\rho>1$ it develops a (macroscopic) shock. Since our results concern the
shock regime, from now on we assume $\lambda\rho>1$.

The microscopic structure of a shock is described by the so called second-class
particle.  Our main result is the explicit expression \eqref{2class-var} for the
variance of the second class particle at any given time $t$. The method,
inspired in work by one of the authors and Fontes \cite{f,ff} for the exclusion
process, relates the mean and the variance of the flux of particles through the
origin with the mean and the variance of the position of a single second class
particle by a linear equation. In this way, if we calculate one, we get the
other. However those works show that for the exclusion process the variance of a
second class particle is asymptotically linear in time; the constant can also be
computed. In the HAD process context, we show how to refine this method to get
an exact formula for the variance of a second class particle at any given time
$t$.

\section{Preliminaries, coupling and Burke's theorem}
\label{preliminaries}

\paragraph{\bf Notation} For a generic point process $M$ we denote
\begin{equation}
  \label{z1}
  \int_0^x f(y,M) M(dy):= \sum_{y\in M\cap[0,x]} f(x,M)
\end{equation}
and recall the Slivnyak-Mecke Theorem: if $M$ is a Poisson process with
intensity $\lambda$, then
\begin{equation}
  \label{z2}
  \bE \Big(\int_0^x f(y,M) M(dy)\Big)= \int_0^x \bE f(y,M) \lambda(y) dy
\end{equation}

\paragraph{\bf The stationary regime and the Burke's theorem} Let $\gamma>0$ and
$S_\chi$ and $W_\chi$ be Poisson processes with intensity $\gamma$ and
$\gamma^{-1}$, respectively. Then $\chi={\mathcal H}(S_\chi,W_\chi,P)$ is a
stationary HAD process with parameter~$\gamma$. Cator and Groeneboom \cite{cg}
proved the following Burke's theorem: $N_\chi$ and $E_\chi$ are independent
Poisson processes with the same density as $S_\chi$ and $W_\chi$, respectively.

\paragraph{\bf Coupling two stationary process} Let $\lambda\rho>1$. Let
$S_\eta$, $W_\sigma$, $I$ and $J$ be mutually independent one-dimensional
Poisson point processes with intensities $\frac{1}{\rho}$, $\frac{1}{\lambda}$,
$\lambda-\frac1\rho$ and $\rho-\frac1\lambda$ respectively. Set
\begin{equation}\label{flux-0}
S_\sigma=S_\eta + I\mbox{ and }W_\eta=W_\sigma + J\,.
\end{equation}
Thus, $S_\sigma$ and $W_\eta$ are two (independent) one-dimensional Poisson
processes with intensity $\lambda$ and $\rho$ respectively. Let $P$ be the
rate-1 two-dimensional Poisson process mentioned in the Introduction.

We run three HAD processes simultaneously with different boundary
conditions but with the same $P$. This is a realization $(\phi,\sigma,\eta)$ of
the so called basic coupling defined by
\[
\phi=\cH(S_\sigma,W_\eta,P)\,,\,\sigma=\cH(S_\sigma,W_\sigma,P)\mbox{ and
}\eta=\cH(S_\eta,W_\eta,P)\,.
\]
For $\chi\in\{\phi,\sigma,\eta\}$ we denote by $\chi(x,t)$ the number of
$\chi$-particles (particles that count for the process $\chi$) in the interval
$[0,x]$ at time $t$, with the convention that particles escaping through a sink
during the time interval $[0,t]$ are located at $0$. The flux of the
discrepancies between $\sigma$ and $\eta$ through the space-time line
$(0,0)$-$(x,t)$ is the process defined by
\begin{equation}
  \label{b9}
  \xi(x,t)=\sigma(x,t)-\eta(x,t)\,.
\end{equation}

There are two kinds of discrepancies: i) those starting at a point $y$ at time
$0$ in $I$ (in the $x$-axis); ii) those starting at a added sink in a time $s$ in $J$ (in
the $t$-axis). The position at time $t$ of the second class particle starting at time zero
at site $y$ is denoted $Z_{(y,0)}(t)$, or shortly $Z_{y}(t)$, and the position
at time $t$ of the second class particle starting at site $0$ at time
$s$ is denoted $Z_{(0,s)}(t)$, or shortly $Z_s(t)$.

The flux of second class particles defined in \reff{b9} is also given by
\begin{eqnarray}
  \label{a44}
  \xi(x,t)&=&\#\big\{y\in I\,:\,Z_{(y,0)}(t)\leq x\big\}
  -\#\big\{s\in J\,:\,Z_{(0,s)}(t)>x\big\}\,.\nonumber\\
  &:=& \xi_+(x,t)-\xi_-(x,t)
\end{eqnarray}

\section{The mean and the variance of  the flux of second class particles}
\label{mean&variance}

In this section we show
the following.
\begin{proposition}\label{mean&var}
For all $x,t\geq 0$
\begin{equation}\label{xi-mean}
\bE\big(\xi(x,t)\big)
=\Bigl(\lambda-\frac{1}{\rho}\Bigr)x - \Bigl(\rho-\frac1\lambda\Bigr)t\,.
\end{equation}
\begin{equation}\label{xi-var}
\var \big(\xi(x,t)\big)=
\Bigl(\lambda-\frac{1}{\rho}\Bigr)x + \Bigl(\rho-\frac1\lambda\Bigr)t\,.
\end{equation}
\end{proposition}

Denote by $\bar M:=|M|$, the total number of points of a generic finite point
process $M$ and $\bchi= \chi(x,t)$. With this notation,
\begin{equation}\label{flux-1}
  \bchi= \bN_\chi+\bW_\chi=\bS_\chi+\bEE_\chi\,,
\end{equation}
for $\chi\in\{\sigma,\eta\}$ and
\begin{equation}
  \label{b10}
  \bxi = \bsigma-\beeta\,.
\end{equation}

\proofof{Proposition \ref{mean&var}} By definition \reff{b9} and Burke's theorem,
\[
\bE\big(\xi(x,t)\big)
=\bE\big(\bxi\big)
=\bE\big(\bsigma\big)-\bE\big(\beeta\big)
=\bE\big(\bN_\sigma+\bW_\sigma-
(\bN_\eta+\bW_\eta)\big)
=\big(\lambda x +\frac{t}{\lambda}\big)-\big(\frac{x}{\rho} + \rho t\big)\,.
\]
This shows \eqref{xi-mean}.

To compute the variance combine first \reff{flux-0} and Burke's theorem to get
\begin{eqnarray}
  \label{b4}
  \cov \big(\bS_\sigma,\bW_\eta\big)&=&\cov \big(\bEE_\sigma,\bN_\eta\big)=0\,.
\end{eqnarray}
Hence, by \reff{flux-1},
\begin{eqnarray}
  \label{b5}
  \cov \big(\bsigma,\beeta\big)&=&\cov \big(\bS_\sigma+\bEE_\sigma,\bN_\eta+\bW_\eta\big)\nonumber\\
&=&\cov \big(\bS_\sigma,\bN_\eta\big)+\cov \big(\bEE_\sigma,\bW_\eta\big)\nonumber\\
&=&\cov \big(\bS_\eta+\bI,\bN_\eta\big)+\cov\big(\bEE_\sigma,\bW_\sigma+\bJ\big)\,\nonumber\\
&=&\cov \big(\bS_\eta,\bN_\eta\big)+\cov \big(\bEE_\sigma,\bW_\sigma\big)\,.\label{var-1}
\end{eqnarray}
Analogous (and simpler) reasoning for $\chi\in\{\eta,\sigma\}$ gives
\begin{equation}
  \label{b13}
  \var \big(\bchi\big)=\cov \big(\bEE_\chi,\bW_\chi\big)+\cov
\big(\bN_\chi,\bS_\chi\big)\,.
\end{equation}

On the other hand, since $\bchi=\bN_\chi+\bW_\chi$,
\begin{equation}
  \label{b7}
 \var \big(\bchi\big)=\var \big(\bN_\chi \big)+\var \big(\bW_\chi \big) + 2\cov \big(\bN_\chi,\bW_\chi\big)\,
\end{equation}
Using $\bW_\chi=\bS_\chi+\bEE_\chi-\bN_\chi$ and $E_\chi$ independent of
$N_\chi$,
$$
\cov \big(\bN_\chi,\bW_\chi\big)=\cov \big(\bN_\chi,\bS_\chi+\bEE_\chi-\bN_\chi)=\cov \big(\bN_\chi,\bS_\chi)-\var \big(\bN_\chi\big)\,.
$$
which implies
\begin{equation}\label{var-2}
\var (\bchi)=\var \big(\bW_\chi\big)-\var \big(\bN_\chi\big) + 2\cov \big(\bN_\chi,\bS_\chi\big)\,.
\end{equation}
Identities \reff{b13} and \reff{var-2} imply
\begin{equation}
  \label{b14}
  \cov \big(\bEE_\chi,\bW_\chi\big)-\cov
  \big(\bN_\chi,\bS_\chi\big) = \var \big(\bW_\chi\big)-\var \big(\bN_\chi\big)
\end{equation}

Finally we compute $\var(\bxi)$ using \reff{var-1} and \reff{b13}:
\begin{eqnarray}
  \label{b15}
  \var \big(\bxi\big)&=&\var \big(\bsigma\big)+\var \big(\beeta\big)-2
  \cov \big(\bsigma,\beeta\big)\nonumber\\
  &=& \cov \big(\bEE_\eta,\bW_\eta\big)-\cov
  \big(\bN_\eta,\bS_\eta\big) + \cov
  \big(\bN_\sigma,\bS_\sigma\big) - \cov \big(\bEE_\sigma,\bW_\sigma\big)\nonumber\\
  &=& \Big(\var \big(\bW_\eta\big)-\var \big(\bN_\eta
  \big)\Big)-\Big(\var \big(\bW_\sigma\big)-\var \big(\bN_\sigma\big)\Big)\nonumber\\
  &=&
  \big(\rho t-\frac{x}{\rho}\big)-\big(\frac{t}{\lambda}-\lambda
  x\big)\,,
\end{eqnarray}
where in the third line we used \reff{b14}. \endproof

\section{The mean and the variance of a second class particle}
\begin{lemma}
  \label{xt+}
For al $x,t\ge 0$,
\begin{eqnarray}
  \label{b1}
  \bE\xi(x,t)_+
  &=&\big(\lambda-\frac{1}{\rho}\big) \Big\{x-\int_0^x\bP\big(Z(t)>z\big)dz\Big\}\,.\label{b31} \\
  \bE\xi(x,t)_-
  &=&  \big(\rho-\frac{1}{\lambda}\big) \Big\{t-\int_0^t\bP\big(Z(u)\leq x\big)du\Big\} \label{b30}
\end{eqnarray}
\end{lemma}
\proofof{Lemma \ref{xt+}} By definition \reff{a44},
\[
\xi(x,t)_+
\,=\,\int_0^x\1\{Z_y(t)\leq x\}{I}(dy)
\]
where $Z_y(t)$ denotes the position at time $t$ of a second-class particle that
starts from the discrepancy at $y\in[0,x]$. By translation invariance and
\eqref{z2},
\begin{eqnarray}
  \label{b1}
  \bE\xi(x,t)_+
  &=&\big(\lambda-\frac{1}{\rho}\big)\int_0^x\bP\big(Z(t)\leq x-y\big)dy\,.\label{den-1}\\
  &=&\big(\lambda-\frac{1}{\rho}\big)\int_0^x\bP\big(Z(t)\leq
  z\big)dz
  \nonumber
\end{eqnarray}
which proves \reff{b31}. Analogously, \reff{b30} follows from
\begin{eqnarray}
  \label{a45}
  \bE\xi(x,t)_- &=& \bE\Bigl(\int_0^t \bP(Z_s(t)>x) J(ds)\Bigr)\,.
\end{eqnarray}
\endproof

The main step to calculate the mean value of $Z$ is the following.
\begin{proposition}\label{normal-1}
For all $x,t\geq 0$
\begin{equation}
  \label{a47}
\int_0^x\bP\big(Z(t)>z\big)dz\,=\,\frac{\rho}{\lambda}\int_0^t\bP\big(Z(u)\leq
x\big)du\,.
\end{equation}
\end{proposition}
\proofof{Proposition \ref{normal-1}} Combining \reff{b31}, \reff{b30} with
\reff{xi-mean},
\begin{eqnarray}
  \label{a46}
  &&\Bigl(\lambda-\frac{1}{\rho}\Bigr)x - \Bigl(\rho-\frac1\lambda\Bigr)t
  \;=\;\bE\big(\xi(x,t)\big)
  \;=\;\bE\big(\xi(x,t)_+\big)-\bE\big(\xi(x,t)_-\big)\nonumber\\
  &&\qquad=\;\big(\lambda-\frac{1}{\rho}\big)\Big\{x-\int_0^x\bP\big(Z(t)> z\big)dz\Big\}\,-\, \big(\rho-\frac{1}{\lambda}\big)\Big\{t-\int_0^t\bP\big(Z(u) \leq
  x\big)du\Big\}\,;\nonumber
\end{eqnarray}
that is,
\begin{equation}
  \label{a48}
  \big(\lambda-\frac{1}{\rho}\big)\int_0^x\bP\big(Z(t)> z\big)dz\;=\; \big(\rho-\frac{1}{\lambda}\big)\int_0^t\bP\big(Z(u) \leq
  x\big)du
\end{equation}
which is equivalent to \reff{a47} for $\lambda\rho> 1$.
\endproof

\proofof{\reff{2class-mean}}
It follows directly from Proposition \ref{normal-1}:
\[
\bE\big(Z(t)\big)
=\lim_{x\to\infty}\int_0^x\bP\big(Z(t)>z\big)dz
=\lim_{x\to\infty}\frac{\rho}{\lambda}\int_0^t\bP\big(Z(u)\,\,\leq x\big)du
=\frac{\rho}{\lambda}t\,.
\]
\endproof

A similar idea works to compute $\var \big(Z(t)\big)$:
\begin{proposition}\label{normal-2}
For all $x,t\geq 0$
\begin{eqnarray}
  \label{f8}
  &&\big(\rho-\frac{1}{\lambda}\big)\Big(t+\int_0^t\bP\big(Z(u)\le x\big)du\Big)
  \\
  && \qquad\qquad =\; \big(\lambda-\frac{1}{\rho}\big)^2\Big(2\int_0^x z
  \bP\big(Z(t)> z\big)dz-\Big(\int_0^x\bP\big(Z(t)>z\big)dz\Big)^2\Big)+\, \var(\xi(x,t)_-)\nonumber
\end{eqnarray}
\end{proposition}

\proofof{\reff{2class-var}} Taking $x\to\infty$ in both members of \reff{f8},
\begin{eqnarray}
  \label{f9}
  \big(\rho-\frac{1}{\lambda}\big)2t
  &=& \big(\lambda-\frac{1}{\rho}\big)^2\var(Z(t))\,,\nonumber
\end{eqnarray}
from where \reff{2class-var} follows. (Notice that, for fixed $t\geq 0$, $\xi_-(x,t)$ decreases to $0$ when $x\to\infty$.)  \endproof

\proofof{Proposition \ref{normal-2}} Noticing that if $\bar{y}\leq y$ and
$Z_y(t)\leq x$ then $Z_{\bar{y}}(t)\leq x$,
\begin{equation}
  \label{a53}
  \xi(x,t)_+^2=\xi(x,t)_+ + 2\int_0^x\int_0^y\1\{Z_y(t)\leq
x\}1\{y\neq \bar y\}I(d\bar{y})I(dy)\,
\end{equation}
as we do for the square of a sum of functions with values in $\{0,1\}$. The
expected value of the double integral in the above equation is given by
\begin{eqnarray}
  \label{a20}
  && \hskip-2cm 2\bE\Big(\int_0^x\int_0^y\1\{Z_y(t)\leq
  x\}1\{y\neq \bar y\}I(d\bar{y})I(dy)\Big)\nonumber\\
  &=&2\big(\lambda-\frac{1}{\rho}\big)^2\int_0^x\int_0^y\bP\big(Z(t)\leq
  x-y\big)d\bar{y}dy\nonumber\\
\label{den-2}
&=&2\big(\lambda-\frac{1}{\rho}\big)^2\int_0^x y\bP\big(Z(t)\leq
x-y\big)dy\,,\nonumber\\
&=&2\big(\lambda-\frac{1}{\rho}\big)^2\Big(x\int_0^x
\bP\big(Z(t)\le z\big)dz -\int_0^x z\bP \big(Z(t)\le z\big)dz \Big)\nonumber\\
&=&2\big(\lambda-\frac{1}{\rho}\big)^2\Big(x^2-x\int_0^x
\bP\big(Z(t)> z\big)dz -\frac{x^2}{2}+\int_0^x z
\bP\big(Z(t)> z\big)dz \Big)\end{eqnarray}
using translation invariance. Take the expectation of \reff{a53} using
\reff{b31} and \reff{den-2} to get
\begin{eqnarray}
  \label{f5}
  \bE\big(\xi(x,t)_+^2\big)
  &=&\big(\lambda-\frac{1}{\rho}\big)\Big(x-\int_0^x\bP\big(Z(t)>z\big)dz\Big)\\
  &&\quad+\,\,\big(\lambda-\frac{1}{\rho}\big)^2\Big(x^2-2x\int_0^x
  \bP\big(Z(t)> z\big)dz +2\int_0^x z
  \bP\big(Z(t)> z\big)dz \Big)\nonumber
\end{eqnarray}
On the other hand, taking the square of \reff{b31},
\begin{equation}
  \label{b12}
  \big(\bE\xi(x,t)_+\big)^2 =\big(\lambda-\frac{1}{\rho}\big)^2
  \Big\{x^2-2x\int_0^x\bP\big(Z(t)> z\big)dz
  +\Big(\int_0^x\bP\big(Z(t)>z\big)dz\Big)^2\Big\}
\end{equation}
Substracting \reff{b12} from \reff{f5}:
\begin{eqnarray}
  \label{f6}
  \var(\xi(x,t)_+)
  &=&  \big(\lambda-\frac{1}{\rho}\big)\Big(x-\int_0^x\bP\big(Z(t)>z\big)dz\Big)
  \nonumber\\
  && \qquad +\, \big(\lambda-\frac{1}{\rho}\big)^2\Big(2\int_0^x z
  \bP\big(Z(t)> z\big)dz-\Big(\int_0^x\bP\big(Z(t)>z\big)dz\Big)^2\Big)\nonumber
\end{eqnarray}
that is, using \reff{xi-var} and  $\var(\xi(x,t))= \var(\xi(x,t)_+) +\var(\xi(x,t)_-)$,
\begin{eqnarray}
  \label{f7}
  \big(\rho-\frac{1}{\lambda}\big)t
  &=&  -\big(\lambda-\frac{1}{\rho}\big)\Big(\int_0^x\bP\big(Z(t)>z\big)dz\Big)
  \nonumber\\
  && \qquad +\, \big(\lambda-\frac{1}{\rho}\big)^2\Big(2\int_0^x z
  \bP\big(Z(t)> z\big)dz-\Big(\int_0^x\bP\big(Z(t)>z\big)dz\Big)^2\Big)\nonumber\\
&& \qquad +\, \var(\xi(x,t)_-)\nonumber
\end{eqnarray}
Use \reff{a47} to get
\begin{eqnarray}
  \label{f7}
  &&\big(\rho-\frac{1}{\lambda}\big)\Big(t+\int_0^t\bP\big(Z(u)\le x\big)du\Big)
  \nonumber\\
  && \qquad =\; \big(\lambda-\frac{1}{\rho}\big)^2\Big(2\int_0^x z
  \bP\big(Z(t)> z\big)dz-\Big(\int_0^x\bP\big(Z(t)>z\big)dz\Big)^2\Big)+\, \var(\xi(x,t)_-)\nonumber
\end{eqnarray}
which shows \reff{f8}. \endproof

\section{The dependence on the initial condition}
Recall that
$$
N(t)=\big(\lambda-\frac{1}{\rho}\big)^{-1}\Big[c_0 t-\Big(\int_0^{tR(S) }S(dz) +\int_0^{tR(W)}W(ds)\Big)\Big]\,
$$
where
$$
c_0=3\Big(\rho-\frac{1}{\lambda}\Big),\, \,R(S)=\frac{1}{\lambda}\Big(\rho-\frac{1}{\lambda}\Big)\,\mbox{ and }R(W)=\frac{1}{\rho}\Big(\rho-\frac{1}{\lambda}\Big)\,.
$$
Thus,
\begin{equation}\label{dep-1}
\bE \big(N(t)\big) = \bE\big(Z(t)\big) = \frac{\rho}{\lambda}t\,\,\,\,\mbox{ and }\,\,\,\,\var \big(N(t)\big)=\var \big(Z(t)\big)=Dt\,
\end{equation}
where
$$
D=\frac{2\big(\rho-\frac{1}{\lambda}\big)}{\big(\lambda-\frac{1}{\rho}\big)^{2}}\,.
$$
(by calculating $N(t$) and using Theorem \ref{2class-mean&var}). The next step is to prove:
\begin{lemma}\label{cov}
$$
\lim_{t\to\infty}\frac{\cov \big(Z(t),N(t)\big)}{t}=D$$
\end{lemma}

\proofof{Lemma \ref{cov}} The proof follows in the same lines of an analogous result
proved by Ferrari \cite{ff} in the TASEP context. To sketch the idea,
denote by $Z^{z|1}(t)$ $\big(\mbox{respectively,} \ Z^{z|0}(t)\big)$ the
position at time $t$ of a second class particle with respect to an initial
configuration conditioned to have a particle at position $z$
(\mbox{respectively,} \ conditioned not to have a particle at position
$z$). Then we can get (by coupling)
\begin{equation} \label{cov1}
\cov \big(Z(t),\int_0^{tR(S) }S(dz)\big) = -\lambda \int_0^{tR(S) } \bE
\big(Z^{z|0}(t) - Z^{z|1}(t)\big)dz
\end{equation}
and
\begin{equation} \label{cov2}
\cov \big(Z(t),\int_0^{tR(W)}W(ds)\big) = -\rho \int_0^{tR(W)} \bE
\big(Z^{s|0}(t) - Z^{'s|1}(t)\big)ds \,.
\end{equation}

Now, we claim that for all $\epsilon > 0$
\begin{equation} \label{pf93} \lim_{t \rightarrow \infty} \sup_{z \in [0,t(R(.)
    - \epsilon )]} |\bE \big(Z^{z|0}(t) - Z^{z|1}(t)\big) -
  \frac{1}{\big(\lambda - \frac{1}{\rho}\big)} |=0 \,.
\end{equation}
To sketch the idea, call (see Section \ref{mean&variance}) $Z_n^{z|1}(t)$ the position at time $t$ of the
second class particles starting at $\xi$ conditioned to have one particle at $z$.
Let $Y_t$ be the position at time $t$ where the processes with initial
configuration being different only at $z$ differ by time $t$. Let $\tau_1$ be the first time the trajectories of $Y_t$ and $Z^{z|1}(t)$ meet each other. Note that $\tau_1$ is finite with probability
one. Before $\tau_1, Z^{z|1}(t) = Z^{z|0}(t)$. After
$\tau_1,Z^{z|1}(t)=Z^{z|1}_1(t)$ and $Z^{z|0}_t=Z^{z|1}_{2}(t)$. Since
$$
\bE
\big(Z^{z|1}_2(t) - Z^{z|1}_{1}(t)\big)=\frac{1}{\big(\lambda -
  \frac{1}{\rho}\big)} \,,
$$
we get \reff{pf93}. Thus, Lemma \ref{cov} follows by combining \reff{cov1},
\reff{cov2} together with \reff{pf93}.
\endproof

\proofof{Theorem \ref{dependence}}
By \reff{dep-1},
$$
\bE\Big[\big(Z(t)-N(t)\big)^2\Big]=2Dt-2\cov \big(Z(t),N(t)\big)\,,
$$
and hence, Theorem \ref{dependence} follows from Lemma \ref{cov}. \endproof

\end{document}